\documentclass[reqno,12pt]{amsart}
\usepackage{amsfonts}
\usepackage{amssymb}
\usepackage{graphicx}
\usepackage{fullpage}
\usepackage{url}

\begin{document}

\newtheorem{theorem}{Theorem}
\newtheorem{lemma}{Lemma}
\newtheorem{proposition}{Proposition}


\def\eqref#1{(\ref{#1})}
\def\eqrefs#1#2{(\ref{#1}) and~(\ref{#2})}
\def\eqsref#1#2{(\ref{#1}) to~(\ref{#2})}
\def\sysref#1#2{(\ref{#1})--(\ref{#2})}

\def\Eqref#1{Eq.~(\ref{#1})}
\def\Eqrefs#1#2{Eqs.~(\ref{#1}) and~(\ref{#2})}
\def\Eqsref#1#2{Eqs.~(\ref{#1}) to~(\ref{#2})}
\def\Sysref#1#2{Eqs. (\ref{#1}),~(\ref{#2})}

\def\secref#1{Sec.~\ref{#1}}
\def\secrefs#1#2{Sec.~\ref{#1} and~\ref{#2}}

\def\appref#1{Appendix~\ref{#1}}

\def\Ref#1{Ref.~\cite{#1}}
\def\Refs#1{Refs.~\cite{#1}}

\def\Cite#1{${\mathstrut}^{\cite{#1}}$}

\def\EQ{\begin{equation}}
\def\endEQ{\end{equation}}


\def\fewquad{\qquad\qquad}
\def\severalquad{\qquad\fewquad}
\def\manyquad{\qquad\severalquad}
\def\manymanyquad{\manyquad\manyquad}

\def\downindex#1{{}_{#1}}
\def\upindex#1{{}^{#1}}

\def\mathtext#1{\hbox{\rm{#1}}}

\def\hp#1{\hphantom{#1}}

\def\parder#1#2{\partial{#1}/\partial{#2}}
\def\parderop#1{\partial/\partial{#1}}

\def\xder#1{{#1}\downindex{x}}
\def\vder#1{{#1}\downindex{v}}

\def\D#1{D\downindex{#1}}

\def\X#1{{\bf X}_{\rm #1.}}
\def\prX#1{{\bf X}_{\rm #1.}^{(1)}}
\def\S#1{S_{\rm #1}}

\def\sgn{{\rm sgn}}

\def\ie/{i.e.}
\def\eg/{e.g.}
\def\etc/{etc.}
\def\const{{\rm const.}}


\title{Symmetry analysis and exact solutions of semilinear heat flow in multi-dimensions}

\author{Stephen C. Anco}
\address{
Department of Mathematics,
Brock University,
St. Catharines, ON L2S 3A1 Canada }
\email{sanco@brocku.ca}

\author{S. Ali}
\address{
School of Electrical Engineering and Computer Sciences,
National University of Sciences and Technology,
H-12 Campus, Islamabad 44000, Pakistan}
\email{sajid\_ali@mail.com}

\author{Thomas Wolf}
\address{
Department of Mathematics,
Brock University,
St. Catharines, ON L2S 3A1 Canada }
\email{twolf@brocku.ca}

\begin{abstract}
A symmetry group method is used to obtain exact solutions for
a semilinear radial heat equation in $n>1$ dimensions
with a general power nonlinearity.
The method involves an ansatz technique to solve
an equivalent first-order PDE system of similarity variables
given by group foliations of this heat equation,
using its admitted group of scaling symmetries.
This technique yields explicit similarity solutions as well as
other explicit solutions of a more general (non-similarity) form
having interesting analytical behavior connected with blow up and dispersion.
In contrast,
standard similarity reduction of this heat equation gives
a semilinear ODE that cannot be explicitly solved by familiar
integration techniques such as point symmetry reduction or integrating factors.
\end{abstract}

\def\sep{;}
\keywords{
semilinear heat equation, similarity reduction, exact solutions,
group foliation, symmetry
}
\subjclass[2010]{35K58\sep 35C06\sep 35A25\sep 58J70\sep 34C14}

\maketitle

\section{ Introduction }

In the study of nonlinear partial differential equations (PDEs),
similarity solutions are important for the understanding of
asymptotic behaviour and attractors, critical dynamics, and blow-up behaviour.
Such solutions are characterized by a scaling homogeneous form
arising from invariance of a PDE under a point symmetry group of
scaling transformations that act on the independent and dependent variables
in the PDE \cite{Olver-book,BlumanAnco-book}.

For scaling invariant PDEs that have only two independent variables,
similarity solutions satisfy an ordinary differential equation (ODE)
formulated in terms of the invariants of the scaling transformations.
However, this ODE can often be very difficult to solve explicitly,
and as a consequence, special ansatzes or ad hoc techniques may be necessary
in order to obtain any solutions in an explicit form. The same difficulties
occur more generally in trying to find explicit group-invariant solutions
to nonlinear PDEs with other types of point symmetry groups.

An interesting example is the semilinear radial heat equation
\EQ\label{heatfloweq}
u_{t} = u_{rr} +(n-1)r^{-1}u_{r} + k |u|^q u,
\quad
k=\pm 1
\endEQ
for $u(t,r)$, which has the scaling symmetry group
\EQ\label{scaling}
t\rightarrow \lambda^2 t,\quad
r\rightarrow \lambda r,\quad
u\rightarrow \lambda^{-2/q}u,
\endEQ
where $r$ denotes the radial coordinate in $n>1$ dimensions.
This equation describes radial heat flow with a nonlinear heat source/sink term
depending on a power $q\neq 0$.
The coefficient $k$ of this term determines the stability of solutions
to the initial-value problem.
In particular,
for $k=-1$ all smooth solutions $u(t,r)$
asymptotically approach a similarity form $u=t^{-1/q}U(r/\sqrt{t})$
exhibiting global dispersive behaviour $u\rightarrow 0$
as $t\rightarrow \infty$ for all $r\geq 0$,
while for $k=1$ some solutions $u(t,r)$
exhibit a blow-up behaviour $u\rightarrow \infty$
given by a similarity form $u=(T-t)^{-1/q}U(r/\sqrt{T-t})$
as $t\rightarrow T <\infty$
\cite{Velazquez,Galaktionov-book}.
In both cases, $U$ satisfies a nonlinear ODE
\EQ\label{similarityeq}
U'' + ((n-1)\xi^{-1} -\tfrac{1}{2}k\xi) U' -\tfrac{1}{2}qk U +k U |U|^{q} =0
\endEQ
with
\EQ
\xi =\begin{cases}
r/\sqrt{t} ,& k=-1 \\
r/\sqrt{T-t} ,& k=1
\end{cases}
\endEQ
which arises from the scaling invariance \eqref{scaling}.
This ODE cannot be explicitly solved by standard integration techniques
\cite{BlumanAnco-book}
such as symmetry reduction or integrating factors when $n\neq1$ and $q\neq0$.
(More specifically,
ODE \eqref{similarityeq} has no point symmetries
$X=\eta(\xi,U)\partial/\partial U +\zeta(\xi,U)\partial/\partial\xi$
and no quadratic first integrals $\Psi=A(\xi,U)+B(\xi,U)U'+C(\xi,U)U'^2$,
as established by solving the standard determining equations
\cite{Olver-book,BlumanAnco-book}
for $X$ and $\Psi$.)
As a result,
few exact solutions $U(\xi)$
other than the explicit constant solution $U=(q/2)^{1/q}$
are apparently known to-date.

In this paper
we will obtain explicit exact solutions for the heat equation \eqref{heatfloweq}
by applying an alternative similarity method
developed in previous work \cite{AncoLiu}
on finding exact solutions to a semilinear radial wave equation
with a power nonlinearity.
The method uses the group foliation equations \cite{Ovsiannikov}
associated with a given point symmetry of a nonlinear PDE.
These equations consist of an equivalent first-order PDE system
whose independent and dependent variables are respectively given by
the invariants and differential invariants of the point symmetry transformation.
In the case of a PDE with power nonlinearities,
the form of the resulting group-foliation system
allows explicit solutions to be found by a systematic separation technique
in terms of the group-invariant variables.
Each solution of the system geometrically corresponds to
an explicit one-parameter family of exact solutions of 
the original nonlinear PDE, such that the family is closed under 
the given symmetry group acting in the solution space of the PDE.

In \secref{method},
we set up the group foliation system
given by the scaling symmetry \eqref{scaling}
for the heat equation \eqref{heatfloweq}
and explain the separation technique that we use 
to find explicit solutions of this system.
The resulting exact solutions of the heat equation are summarized
in \secref{results}.
These solutions include explicit similarity solutions
as well as other solutions whose form is not scaling homogeneous,
and we discuss their analytical features of interest
pertaining to blow-up and dispersion.
Finally, we make some concluding remarks in \secref{remarks}.

Related work using a similar method applied to nonlinear diffusion equations
appears in \Ref{QuZhang}. 
Group foliation equations were first used successfully 
in \Refs{Golovin,NutkuSheftel,Sheftel} 
for obtaining exact solutions to nonlinear PDEs 
by a different method that is applicable 
when the group of point symmetries of a given PDE is infinite-dimensional,
compared to the example of a finite-dimensional symmetry group 
considered both in \Ref{AncoLiu} and in the present work.

\section{ Symmetries and group foliation }
\label{method}

For the purpose of symmetry analysis and finding exact solutions,
it is easier to work with a slightly modified form of the heat equation \eqref{heatfloweq}:
\EQ\label{heateq}
u_{t} = u_{rr} +(n-1)r^{-1}u_{r} + k u^{q+1},
\quad
k=\const,\quad
q\neq 0.
\endEQ
In $n>1$ dimensions,
this heat equation \eqref{heateq} admits only two point symmetries:
\begin{align}
&\mathtext{ time translation } \quad
\X{trans} =\parderop{t} \quad\mathtext{ for all $q$},
\label{transsymm}\\
&\mathtext{ scaling } \quad
\X{scal} =2t\parderop{t} + r\parderop{r} -(2/q) u\parderop{u}
\quad\mathtext{ for all $q\neq 0$},
\label{scalsymm}
\end{align}
where $\X{}$ is the infinitesimal generator of a one-parameter group of
point transformations acting on $(t,r,u)$.
There are no special powers or dimensions for which any extra point symmetries
exist for equation \eqref{heateq},
as found by a direct analysis of the symmetry determining equations.

To proceed with setting up the group foliation equations
using the scaling point symmetry,
we first write down the invariants and differential invariants
determined by the generator \eqref{scalsymm}.
The simplest invariants in terms of $t,r,u$ are given by
\EQ
x=t/r^2 ,\quad
v=u/r^{p},
\label{scalxv}
\endEQ
satisfying $\X{scal}x=\X{scal}v=0$
with
\EQ
p=-2/q.
\endEQ
A convenient choice of differential invariants satisfying
$\prX{scal}G=\prX{scal}H=0$ for $G(t,r,u_{t})$ and $H(t,r,u_{r})$ consists of
\EQ
G=r^{2-p} u_{t} ,\quad
H=r^{1-p} u_{r} ,
\label{scalGH}
\endEQ
where
$\prX{scal}=\X{scal}-(2+2/q) u_t\parderop{u_t} -(1+2/q) u_r\parderop{u_r}$
is the first-order prolongation of the generator \eqref{scalsymm}.
Here $x$ and $v$ are mutually independent,
while $G$ and $H$ are related by equality of mixed $r,t$ derivatives
on $u_{t}$ and $u_{r}$, which gives
\EQ
\D{r}( r^{p-2} G ) = \D{t}( r^{p-1} H )
\label{GHmixedeq}
\endEQ
where $\D{r},\D{t}$ denote total derivatives with respect to $r,t$.
Furthermore,
$v,G,H$ are related through the heat equation \eqref{heateq} by
\EQ
r^{p-2} G -\D{r}( r^{p-1} H ) = r^{p-2}( (n-1)H + kv^{q+1} ) .
\label{GHheateq}
\endEQ
Now we put $G=G(x,v)$, $H=H(x,v)$ into equations \eqrefs{GHmixedeq}{GHheateq}
and use equation \eqref{scalxv} combined with the chain rule
to arrive at a first-order PDE system
\begin{align}
&(p-2)G -pv\vder{G} -2x\xder{G} -\xder{H} + H\vder{G} -G\vder{H} =0 ,
\label{scalGHeq1}\\
& G -(p+n-2)H +pv\vder{H} +2x\xder{H} -H\vder{H} =kv^{q+1} ,
\label{scalGHeq2}
\end{align}
with independent variables $x,v$, and dependent variables $G,H$.
These PDEs are called the {\it scaling-group resolving system}
for the heat equation \eqref{heateq}.

The respective solution spaces of 
equation \eqref{heateq} and system \sysref{scalGHeq1}{scalGHeq2}
are related by a group-invariant mapping that is defined through 
the invariants \eqref{scalxv} and differential invariants \eqref{scalGH}. 
In particular, 
the map $(G,H) \rightarrow u$ is given by 
integration of a consistent pair of parametric first-order ODEs
\EQ
u_{t} =r^{p-2} G(t/r^2,u/r^{p}) ,\quad
u_{r} =r^{p-1} H(t/r^2,u/r^{p})
\label{umap}
\endEQ
whose general solution will involve a single arbitrary constant. 
The inverse map $u\rightarrow (G,H)$ can be derived 
in the same way as shown in \Ref{AncoLiu} for the wave equation,
which gives the following correspondence result. 

\begin{lemma}
Solutions $(G(x,v),H(x,v))$ of 
the scaling-group resolving system \sysref{scalGHeq1}{scalGHeq2}
are in one-to-one correspondence with 
one-parameter families of solutions $u(t,r,c)$ of 
the heat equation \eqref{heateq}
satisfying the scaling-invariance property
\EQ
\lambda^{-p}u(\lambda^2 t,\lambda r,c) = u(t,r,\tilde c(\lambda,c)) 
\label{orbit}
\endEQ
where $\tilde c(1,c)=c$ in terms of an arbitrary constant $c$. 
\end{lemma}

This correspondence leads to an explicit characterization of 
similarity solutions of the heat equation \eqref{heateq}
in terms of a condition on solutions of 
the scaling-group resolving system \sysref{scalGHeq1}{scalGHeq2}. 
Consider any one-parameter family of solutions
\EQ\label{scalu}
u(t,r)= r^p v,\quad
v = V(x,c) ,
\endEQ
having a scaling-homogeneous form,
where
\EQ\label{scalUeq}
4 x^2 V'' -(1+(2p+n-4)2x) V' +p(p+n-2)V + k V^{q+1} =0
\endEQ
is the ODE given by reduction of PDE \eqref{heateq}.
From relation \eqref{scalGH} we have
\EQ
G(x,V(x,c))=V'(x,c) ,\quad
H(x,V(x,c))=pV(x,c) - 2xV'(x,c) . 
\label{GHUrel}
\endEQ
Next we eliminate $c$ in terms of $x$ and $v$ 
by using the implicit function theorem on $V(x,c)-v=0$ to express
$c=C(x,v)$. 
Substitution of this expression into equation \eqref{GHUrel} yields
\EQ
H+2x G =pv
\label{scalinvGHsol}
\endEQ
where $G=V'(x,C(x,v))$, $H=pv-2xV'(x,C(x,v))$ are some functions of $x,v$. 
The relation \eqref{scalinvGHsol} 
is easily verified to satisfy PDE \eqref{scalGHeq1}.
In addition, PDE \eqref{scalGHeq2} simplifies to
\EQ
-4 x^2( \xder{G} + G\vder{G} ) + (1+(2p+n-4)2x) G  = p(p+n-2) v +k v^{q+1} .
\label{scalinvGHeq}
\endEQ
We then see that the characteristic ODEs for solving this first-order PDE
are precisely
\EQ
dv/dx =G ,\quad
-4x^2 dG/dx + (1+(2p+n-4)2x) G =p(p+n-2) v +k v^{q+1} ,
\endEQ
which are satisfied due to equations \eqrefs{scalUeq}{GHUrel}.
Hence, we have established the following result.

\begin{lemma}
There is a one-to-one correspondence between 
one-parameter families of similarity solutions \eqref{scalu} 
of heat equation \eqref{heateq}
and solutions of the scaling-group resolving system \sysref{scalGHeq1}{scalGHeq2}
that satisfy the similarity relation \eqref{scalinvGHsol}.
\end{lemma}

We now note that, under the mapping \eqref{umap},
static solutions $u(r)$ of the heat equation
correspond to solutions of the scaling-group resolving system with $G=0$.
Consequently, hereafter we will be interested only in solutions
such that $G\neq0$,
corresponding to dynamical solutions of the heat equation.

To find explicit solutions of the PDE system \sysref{scalGHeq1}{scalGHeq2}
for $G(x,v),H(x,v)$,
we will exploit its following general features.
First,
the power nonlinearity $k u^{q+1}$ in the heat equation
appears only as an inhomogeneous term $k v^{q+1}$ in the PDE \eqref{scalGHeq2}.
Second,
in both PDEs \eqrefs{scalGHeq1}{scalGHeq2}
the linear terms that involve $v$ derivatives
have the scaling homogeneous form $vG_v$ and $vH_v$ with respect to $v$.
Third,
the nonlinear terms in the homogeneous PDE \eqref{scalGHeq1}
have the skew-symmetric form $HG_v-GH_v$,
while $HH_v$ is the only nonlinear term appearing
in the non-homogeneous PDE \eqref{scalGHeq2}.
These features suggest that this PDE system can be expected to have solutions
given by the separable power form
\EQ
G=g_1(x) v^{a} + g_2(x) v,\quad
H=h_1(x) v^{a} + h_2(x) v,\quad
a\neq 1 .
\label{GHansatz}
\endEQ
For such an ansatz,
we readily see that the linear derivative terms $G_x$, $H_x$, $vG_v$, $vH_v$
in each PDE \eqrefs{scalGHeq1}{scalGHeq2} will contain the same powers $v,v^a$
that appear in both $G$ and $H$,
and moreover the nonlinear term $HG_v-GH_v$
in the homogeneous PDE \eqref{scalGHeq1} will produce only the power $v^a$
due to the identities $v^a(v)_v- v(v^a)_v = (a-1)v^a$
and $v(v)_v- v(v)_v = v^a(v^a)_v- v^a (v^a)_v = 0$.
Similarly we see that the nonlinear term $HH_v$
in the non-homogeneous PDE \eqref{scalGHeq2} will only yield the powers
$v$, $v^a$, $v^{2a-1}$.
Since we have $a\neq1$ and $q\neq0$,
the inhomogeneous term $k v^{q+1}$ must therefore balance one of the powers
$v^{2a-1}$ or $v^a$.

In the case when we balance $q+1=a$,
the terms containing $v^a=v^{q+1}$ and $v^{2a-1}$ in PDE \eqref{scalGHeq2}
immediately yield
\EQ
h_1=0, \quad
g_1=k .
\endEQ
Then the terms containing $v^a$ in PDE \eqref{scalGHeq1} reduce to
\EQ
h_2=0
\endEQ
which leads to a simplification of the remaining terms in both PDEs,
yielding
\EQ
x g_2'=g_2=0 .
\endEQ
Thus, for this case,
the ansatz \eqref{GHansatz} gives a one-term solution for $G$ with $H=0$.

In the other case, balancing $2a-1=q+1$, we get
\EQ
a=1+\frac{q}{2} .
\endEQ
The terms containing  $v$, $v^a$, $v^{2a-1}=v^{q+1}$
in the PDEs \eqref{scalGHeq1} and \eqref{scalGHeq2}
then yield
\begin{subequations}\label{sys1eqs}
\begin{align}
& h_1' + 2x g_1' +(1-\frac{q}{2}h_2)g_1 + \frac{q}{2} h_1 g_2 =0
\label{eq1}\\
& h_2' + 2x g_2'+ 2 g_2 =0
\label{eq2}\\
& 2x h_1' -(n-1+\frac{q+4}{2}h_2)h_1 + g_1 =0
\label{eq3}\\
& 2x h_2' -h_2{}^2 -(n-2)h_2 + g_2 =0
\label{eq4}\\
& (1+\frac{q}{2}) h_1{}^2+k =0 .
\label{h1eq}
\end{align}
\end{subequations}
Through equations \eqref{h1eq}, \eqref{eq3} and \eqref{eq4},
we obtain
\begin{align}
& h_1=\pm\sqrt{\frac{-2k}{q+2}}
\label{h1}\\
& g_1 =\pm\sqrt{\frac{-2k}{q+2}} (n-1+\frac{q+4}{2}h_2)
\label{g1}\\
& g_2 = -2x h_2' +(n-2 +h_2)h_2
\label{g2}
\end{align}
and thereby we find that equations \eqref{eq1} and \eqref{eq2} reduce to
an overdetermined system of nonlinear ODEs
\begin{subequations}\label{sys1}
\begin{align}
& 4x^2 h_2' - 2x h_2{}^2 -(1+2(n-2)x)h_2 =c=\const
\label{h2ode1}\\
& 4x h_2' -\frac{q(2+q)}{4} h_2{}^2+ 2 h_2 +n-1=0
\label{h2ode2}
\end{align}
\end{subequations}
where $c$ is an arbitrary integration constant.
The ODE system \eqref{sys1} can be solved by
a systematic integrability analysis,
which we have carried out using computer algebra
(discussed in more detail in \secref{integrabilityruns}). 
The results of the analysis give three two-term solutions with $G\neq0$.

\begin{proposition}
For $n\neq1$,
ansatz \eqref{GHansatz} yields altogether four solutions of
the scaling-group resolving system \sysref{scalGHeq1}{scalGHeq2}
with $G\neq 0$:
\begin{align}
&
G = k v^{q+1},\quad
H = 0 , \quad
q\neq -1
\label{scalingsol1}\\
&
G= \pm(4-n) \sqrt{\frac{-k (n-2)} {n-3}} v^ {(n-3)/(n-2)},\quad
H= \frac{1}{4-n} G + (2-n) v,
\label{scalingsol2}\\
&\qquad
q=\frac{2}{2-n}\neq -1 , \quad
n\neq 2,3,4
\nonumber\\
&
G= \frac{3 v}{3 x+1} \pm\frac{3\sqrt{k}}{2 v},\quad
H= \frac{2}{3} G - \frac{v}{2},\quad
q=-4,\quad
n=5/2
\label{scalingsol3}\\
&
G = \frac{3 v \left (1 \pm \sqrt{-2 k} v \right)}{3 x +1},\quad
H = \frac{1}{6} (3 x+1) G -\frac{3 x-1}{3 x+1} v,\quad
q=2,\quad
n=5/2
\label{scalingsol4}
\end{align}
None of these solutions satisfy the similarity relation \eqref{scalinvGHsol}.
\end{proposition}

Motivated by the success of the ansatz \eqref{GHansatz},
we now consider a more general three-term ansatz
\EQ
G=g_1(x) v^a + g_2(x) v^b +g_3(x)v ,\quad
H=h_1(x) v^a + h_2(x) v^b +h_3(x)v
\label{3termGHansatz}
\endEQ
where
\EQ
a\neq b,\quad
a\neq1,\quad
b\neq1 .
\label{GHconds}
\endEQ
This ansatz leads to a more complicated analysis
compared to the previous two-term ansatz.
Specifically,
the homogeneous PDE \eqref{scalGHeq1} now contains the power $v^{a+b-1}$
in addition to $v,v^a,v^b$, while the non-homogeneous PDE \eqref{scalGHeq2}
contains the further powers $v^{2a-1},v^{2b-1},v^{q+1}$.
We determine the exponents in these powers by a systematic examination of
all possible balances.

Firstly,
since $q\neq0$ in PDE \eqref{scalGHeq2},
$v^{q+1}$ must balance one of $v^a,v^{2a-1},v^{a+b-1}$.
(Note, by the symmetry $a\leftrightarrow b$ in the ansatz \eqref{3termGHansatz},
the other possibilities $v^b,v^{2b-1}$ for balancing $v^{q+1}$ are redundant.)
Secondly,
$v^{2a-1}$ can balance only $v^{q+1}$ or $v^b$ due to conditions \eqref{GHconds},
and otherwise if $v^{2a-1}$ is unbalanced then
its coefficient $a h_1{}^2$ must vanish.
Likewise $v^{2b-1}$ can balance only $v^{q+1}$ or $v^a$,
and otherwise its coefficient $b h_2{}^2$ must vanish.
In a similar way,
either $v^{a+b-1}$ balances $v^{q+1}$ or $v$,
and otherwise if $v^{a+b-1}$ is unbalanced then
its coefficient $(a+b)h_1 h_2$ vanishes.
Finally, in PDE \eqref{scalGHeq1},
$v^{a+b-1}$ can balance only $v$,
and otherwise the factor $g_1 h_2-g_2 h_1$
must vanish in the coefficient of $v^{a+b-1}$.

Several cases arise from examining all of these different possibilities.
After eliminating all trivial cases that lead to $h_2=g_2=0$
(whereby the ansatz \eqref{3termGHansatz} just reduces to
the previously considered two-term case \eqref{GHansatz}),
we find the following non-trivial cases to consider:
\begin{align}
&
q=a=2,\quad
b=0;
\label{case1}\\
&
q=-3/2,\quad
a=0,\quad
b=-1/2;
\label{case2}\\
&
q=-2/3,\quad
a=-b=-1/3.
\label{case3}
\end{align}

For case \eqref{case1},
the PDEs \eqref{scalGHeq1} and \eqref{scalGHeq2} yield
\begin{subequations}\label{sys2eqs}
\begin{align}
& h_1' + 2x g_1' +(1-h_3)g_1 + h_1 g_3 =0
\label{2ndeq1}\\
& 2x h_1' -(n-1+3h_3)h_1 + g_1 =0
\label{2ndeq2}\\
& h_3' + 2x g_3' +2 g_2 h_1-2 g_1 h_2 + 2 g_3 =0
\label{2ndeq3}\\
& 2x h_3' -h_3{}^2 -2 h_1 h_2 +(2-n)h_3 + g_3 =0
\label{2ndeq4}\\
& h_2' + 2x g_2' +(3+h_3)g_2 - g_3 h_2 =0
\label{2ndeq5}\\
& 2x h_2' -h_2 h_3 +(3-n)h_2 +g_2 =0
\label{2ndeq6}\\
& 2 h_1{}^2+k =0 .
\label{2ndeq7}
\end{align}
\end{subequations}
From equations \eqref{2ndeq7}, \eqref{2ndeq6}, \eqref{2ndeq4}, \eqref{2ndeq2},
we have
\begin{align}
& h_1=\pm\sqrt{\frac{-k}{2}}
\label{2ndh1}\\
& g_2 = -2x h_2' +(n-3 +h_3)h_2
\label{2ndg2}\\
& g_3 = -2x h_3' +(n-2 +h_3)h_3 \pm\sqrt{-2k} h_2
\label{2ndg3} \\
& g_1 =\pm\sqrt{\frac{-k}{2}} (n-1+3h_3)
\label{2ndg1}
\end{align}
and then equation \eqref{2ndeq1} gives
\EQ
h_2= \pm\frac{1}{\sqrt{-2k}}( 4x h_3' -2h_3{}^2 +2h_3 +n-1 )
\label{2ndh2}
\endEQ
The remaining equations \eqref{2ndeq3} and \eqref{2ndeq5} become, respectively,
\begin{subequations}\label{sys2}
\begin{align}
& 4x^2 h_3'' -(12x h_3 +2(n-6)x +1) h_3' +4 h_3{}^3 -6 h_3{}^2 +2(3-2n)h_3 =0
\label{h3ode1}\\
& 4x^3 h_3''' -x(4x h_3 +2(n-13)x +1) h_3''
+(6x h_3{}^2 +(2(n-12)x +1)h_3 +9(4-n)x -3/2) h_3'
\nonumber\\&\qquad
-12x^2 h_3'{}^2 -\tfrac{1}{2}(2h_3{}^2-2 h_3 +1-n)(h_3{}^2- 2 h_3 +5-2n)=0
\label{h3ode2}
\end{align}
\end{subequations}
which is an overdetermined system of two nonlinear ODEs for $h_3(x)$.
We solve this system \eqref{sys2} by an integrability analysis
using computer algebra.
This yields one solution with $G=0$,
plus two solutions with $G\neq 0$
which are summarized in Proposition~2.

In a similar way, each of the cases \eqref{case2} and \eqref{case3}
leads to an overdetermined system of four nonlinear ODEs for $h_2(x),h_3(x)$.
For both cases the results of an integrability analysis yield
only solutions with $G=0$.

Thus we have the following result.

\begin{proposition}
For $n\neq1$,
ansatz \eqref{3termGHansatz} yields two additional solutions of
the scaling-group resolving system \sysref{scalGHeq1}{scalGHeq2}
with $G\neq0$:
\begin{align}
&
G= \pm\frac{3}{4}\sqrt{-2 k} \left (v \pm \frac{1}{\sqrt{-2 k}} \right )^2,\quad
H= \frac{2}{3}G + v \pm \frac{2}{\sqrt{-2 k}},\quad
q=2, \quad n=5/2
\label{scalingsol5}\\
&
G= \pm \frac{15}{4} \sqrt{-2k} \left (v \mp \frac{1}{\sqrt{-2 k}} \right )^2,\quad
H= \frac{2}{15}G +v \mp \frac{2}{\sqrt{-2 k}},\quad
q=2, \quad n=5/2
\label{scalingsol6}
\end{align}
Neither of these solutions satisfies the similarity relation \eqref{scalinvGHsol}.
\end{proposition}

\subsection{Computational remarks}
\label{integrabilityruns}
\noindent\newline\indent
The integrability analysis of the previous ODE systems is non-trivial
due to the degree of nonlinearity of the ODEs
and the algebraic complexity of the coefficients 
in addition to the appearance of parameters in each system. 

For the first ODE system \eqref{sys1}, 
the integrability analysis consists of the following main steps. 
We eliminate $h_2'$ to get a single algebraic equation, 
which is quadratic in $h_2$.
By differentiating this equation and using it to eliminate $h_2'$ 
from either of the original ODEs in the system, 
we obtain a second algebraic equation, 
which is cubic in $h_2$.
The coefficients in each algebraic equation are expressions in terms of 
the independent variable $x$ and the parameters $q,n$. 
We next use cross-multiplication repeatedly to eliminate 
the highest-degree monomial terms in both of the algebraic equations 
until one equation no longer contains $h_2$ 
while the other equation is linear in $h_2$. 
At each algebraic elimination step, 
we must note that a case distinction will arise if 
the coefficient of a highest-degree monomial vanishes 
for some values of $q$ or $n$. 
For each case, once the final algebraic equations have been obtained, 
we solve the equation without $h_2$ by splitting it with respect to $x$,
which will yield conditions on the parameters $q,n$,
and we then solve the linear equation for $h_2$ 
subject to these conditions (if any). 

The integrability analysis for the second ODE system \eqref{sys2}
is the same except that we must first use differentiation 
combined with cross-multiplication to eliminate $h_3'''$ and $h_3''$,
thus reducing the differential order of the system down to first-order,
where the coefficients are expressions in terms of 
the independent variable $x$ and the single parameter $n$. 
We may then proceed as before by using algebraic elimination 
to reduce this system to a linear equation that can be solved for $h_3$ 
and an equation that does not contain $h_3$ and thereby determines $n$. 

A similar integrability analysis applies to the two other ODE systems, 
each of which requires solving four ODEs that contain two dependent variables 
$h_2,h_3$ and their derivatives up to second order, in addition to 
the independent variable $x$ and the single parameter $n$. 

Because of the complexity of the algebraic expressions 
and the number of case distinctions that arise in these analyses, 
it is very difficult for an automatic computer algebra program
to fully classify and find all solutions.
(For example,
the Maple program {\sc RiffSimp} running on a workstation for several days 
was unable to complete the full computation for 
any of the second, third, and fourth systems.)

To overcome these difficulties,
we have used the interactive package {\sc Crack} \cite{crack}
which has a wide repertoire of techniques available, including
eliminations, substitutions, integrations, length-shortening of equations,
and factorizations, among others. 
Using {\sc Crack}, 
the complete solution of the integrability analysis was obtained 
in about 50 interactive steps taking 2 seconds in total 
for the system \eqref{sys1}, 
and about 400 interactive steps taking 7 minutes in total 
for the system \eqref{sys2},
while less than 100 steps taking under 1 minute in total 
were needed for each of the other two systems \cite{runs}.

\section{ Exact solutions }
\label{results}

To obtain explicit solutions $u(t,r)$ of the heat equation \eqref{heateq}
from solutions $(G(x,v),H(x,v))$
of its scaling-group resolving system \sysref{scalGHeq1}{scalGHeq2},
we integrate the corresponding pair of
parametric first-order ODEs \eqref{scalGH}.
The integration yields a one-parameter solution family $u(t,r,c)$
which is closed under the action of the group of scaling transformations
\eqref{scaling}.

\begin{theorem}
The semilinear heat equation \eqref{heateq} has the following exact solutions
arising from the explicit solutions of its scaling-group resolving system
found in Propositions~1 and~2:
\begin{align}
&
u= (-k q (t+c))^{-1/q} ,
\quad
q\neq 0
\label{usol1}\\
&
u= \left (\pm \sqrt{\frac{-k}{(n-2)(n-3)}}
\left (\frac{r}{2}-\frac{(n-4)(t+c)}{r}\right )\right )^{n-2} ,
\quad
q= \frac{2}{2-n}\neq -1 ,\quad
n\neq 2,3,4
\label{usol2}\\
&
u = \left( \pm\sqrt{k}\left( 1+ c(3 t+r^{2}) \right) \left(\frac{3 t}{r}+r\right) \right)^{1/2} ,
\quad
q=-4, \quad n=5/2
\label{usol3}\\
&
u= \pm\frac{5(3 t+ r^{2})}{\left( r(15 t + r^{2})+ c\sqrt{r}\right)\sqrt{-2 k}},
\quad
q=2, \quad n=5/2
\label{usol4}\\
&
u =\pm \frac{3(t + c -r^{2})}{r(3 (t+c) + r^{2})\sqrt{-2 k}},\quad
q=2, \quad n=5/2
\label{usol5}\\
&
u= \pm\frac{5(3 (t+c) + r^{2} )}{r(15 (t+c) + r^{2})\sqrt{-2 k}} ,
\quad
q=2, \quad n=5/2
\label{usol6}
\end{align}
where $c$ is an arbitrary constant.
\end{theorem}

Modulo time-translations $t\rightarrow t-c$,
solutions \eqref{usol1}, \eqref{usol2}, \eqref{usol5}, \eqref{usol6}
are similarity solutions
since their form with parameter $c=0$ is preserved
under scaling transformations \eqref{scaling} on $r,t,u$.
In contrast,
the parameter $c$ in solutions \eqref{usol3} and \eqref{usol4}
cannot be removed by time-translations,
and consequently the form of these solutions is not scaling-homogeneous
since $c$ gets scaled under the transformations \eqref{scaling} on $r,t,u$.
Thus,
with respect to the action of the full group of
point symmetries generated by time-translations and scalings
for the heat equation \eqref{heateq},
the solutions \eqref{usol3} and \eqref{usol4} for $c\neq 0$ yield
two-dimensional orbits of non-similarity solutions given by
\begin{align}
&
u = \left( \pm\sqrt{k}\left(1+ c(3(t+{\tilde c})+r^{2})\right) \left(\frac{3 (t+{\tilde c})}{r}+r \right) \right)^{1/2} ,
\quad
q=-4, \quad n=5/2
\label{nonsimsol1}\\
&
u= \pm\frac{5(3 (t+{\tilde c})+ r^{2})}{\left(r (15(t+{\tilde c}) + r^{2}) + c\sqrt{r}\right)\sqrt{-2 k}},
\quad
q=2, \quad n=5/2
\label{nonsimsol2}
\end{align}
whereas the solutions \eqref{usol1}, \eqref{usol2}, \eqref{usol5}, \eqref{usol6}
represent one-dimensional time-translation orbits of similarity solutions
\begin{align}
&
u= (-k q)^{-1/q} r^{-2/q} (t/r^2)^{-1/q} ,
\quad
q\neq 0
\label{simsol1}\\
&
u= \left(\frac{(n-2)(n-3)}{-4k}\right)^{1-n/2} r^{n-2} \left(\pm (1-2(n-4)t/r^2 )\right)^{n-2} ,
\quad
q= \frac{2}{2-n}\neq -1 ,\quad
n\neq 2,3,4
\label{simsol2}\\
&
u =\pm(-2 k)^{-1/2} r^{-1} \frac{3(t/r^2-1)}{3t/r^2+1},
\quad
q=2, \quad n=5/2
\label{simsol3}\\
&
u= \pm(-2 k)^{-1/2} r^{-1} \frac{5(3t/r^2+1)}{15t/r^2+1},
\quad
q=2, \quad n=5/2.
\label{simsol4}
\end{align}
Note that solution \eqref{usol6} is a special case of
solution \eqref{nonsimsol2} given by $c=0$.
We remark that solution \eqref{nonsimsol2} previously has been obtained
in work \cite{Sophocleous} on nonlinear diffusion equations 
through use of Bluman and Cole's nonclassical method. 

Among all these solutions,
the ones \eqref{usol1} and \eqref{usol2} that exist for integer values of $n$
describe radial heat flow in $\mathbb{R}^n$.
In section~\secref{integercase} we will discuss their analytical features
related to blow-up and dispersion for $n\neq 1$.
The remaining solutions
\eqref{usol5}, \eqref{usol6}, \eqref{nonsimsol1}, \eqref{nonsimsol2}
that exist only for a non-integer value $n=5/2$
have a different interpretation
describing heat flow in the plane $\mathbb{R}^2$ with a point-source of
radial heat flux at the origin,
as we will show later in section~\secref{nonintegercase}.
We will also discuss this interpretation for the solution \eqref{usol2}
in the case of non-integer values of $n$.

Before proceeding, we observe that the heat equation \eqref{heateq}
for all values of $n$ can be written in the form of a gradient flow
\EQ
u_{t}=-\delta E/\delta u
\label{gradientflow}
\endEQ
using the ``energy'' integral
\EQ
E=\int_0^\infty \left ( \frac{1}{2} u_{r}^{2}- k f(u) \right ) r^{n-1}dr
\label{energy}
\endEQ
with
\EQ
f(u)=\begin{cases}
\dfrac{1}{q+2} u^{q+2},
& q \neq -2,\\
\ln|u|,
& q=-2.
\end{cases}
\endEQ
Here $\delta/\delta u$ denotes the usual variational derivative
(i.e. Euler operator) with respect to $u$.
This integral \eqref{energy} obeys the equation
\EQ
\frac{dE}{dt} =
-\displaystyle\lim_{r\to 0}(r^{n-1}u_{r}u_{t})
-\int_0^\infty u_{t}^{2} r^{n-1}dr
\label{energyflux}
\endEQ
for any solution $u(t,r)$ with sufficient asymptotic decay for large $r$.
If the ``energy flux'' $\displaystyle\lim_{r\to 0} (r^{n-1}u_{r}u_{t})$
of a solution $u(t,r)$ is non-negative then equation \eqref{energyflux}
shows that $E$ is a decreasing function of $t$.
As a result,
for solutions $u(t,r)$ that also are non-negative,
both terms in $E$ will be non-negative if $k<0$ and $q\neq -2$,
so then $E$ decreases to zero as $t\rightarrow\infty$.
In this case $u(t,r)$ will have dispersive behaviour such that
$u\rightarrow 0$ and $u_{r}\rightarrow 0$ for all $r\geq0$
as $t\rightarrow\infty$.
If instead $k>0$ or $q=-2$,
then the two terms in $E$ will have opposite signs,
in which case $E$ may decrease without bound,
allowing $u(t,r)$ to have blow-up behaviour such that
$u\rightarrow\infty$ or $u_{r}\rightarrow\infty$
as $t\rightarrow T$ for some $T<\infty$.

\subsection{ Behavior of solutions for $n=2,3,$ etc. }
\label{integercase}
\noindent\newline\indent
Similarity solution \eqref{usol1} is spatially homogeneous
and has no restriction on the sign of $k$.
It thus represents the general solution of the ODE $u_t=ku^{q+1}$
with an arbitrary nonlinearity power $q+1\neq 1$.
For $q>0$ the behaviour of $u$ is determined by the sign of $k$.
In the case $k<0$,
$u= (|k|q)^{-1/q}/(t+c)^{1/q} \rightarrow 0$ is dispersive
as $t\rightarrow \infty$,
whereas in the case $k>0$,
$u= (kq)^{-1/q}/(T-t)^{1/q} \rightarrow \infty$ has a blow-up
for $t\rightarrow T$ with $0< T=-c<\infty$.

\begin{figure}[h]
\begin{center}
\includegraphics[scale=0.5]{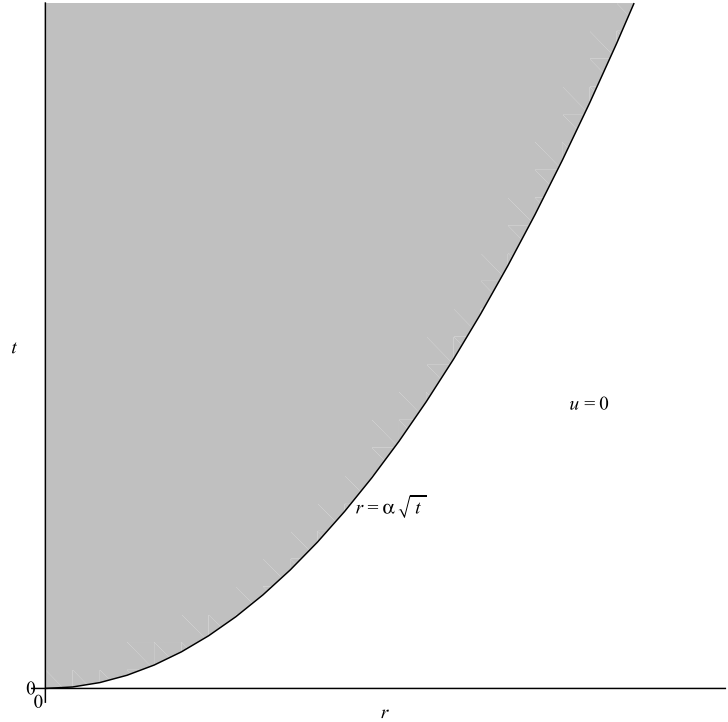}
\end{center}
\caption{Space-time graph for solution \eqref{usol2'}.}\label{eq-sol2-space-time}
\end{figure}

In contrast,
similarity solution \eqref{usol2} is restricted to the special nonlinearity power
$q+1=(n-4)(n-2)\geq 1/3$
and requires $k<0$ and $n\geq 5$.
In this solution, for all $t>0$,
$u$ is singular at $r=0$ and is unbounded as $r\rightarrow \infty$.
Interestingly,
$u$ vanishes on the space-time parabola given by $t+c=r^2/(2n-8) \geq 0$.
As a consequence we can modify the solution to have better analytical behaviour
by using this parabola (with parameter $c=0$) as a cutoff
such that
\EQ
u(t,r) =\begin{cases}
0, & r\geq \alpha \sqrt{t} \\
\beta r^{2-n} \left ( t-(r/\alpha)^2 \right)^{n-2}, & 0\leq r\leq \alpha \sqrt{t}\end{cases}
\label{usol2'}
\endEQ
(see figure \ref{eq-sol2-space-time}) with
\EQ
\alpha = \sqrt{2(n-4)} ,
\quad
\beta =\left(\dfrac{(n-2)(n-3)}{|k|(n-4)^2}\right)^{1-n/2} .
\endEQ

\begin{figure}[h]
\begin{center}
\includegraphics[scale=0.5]{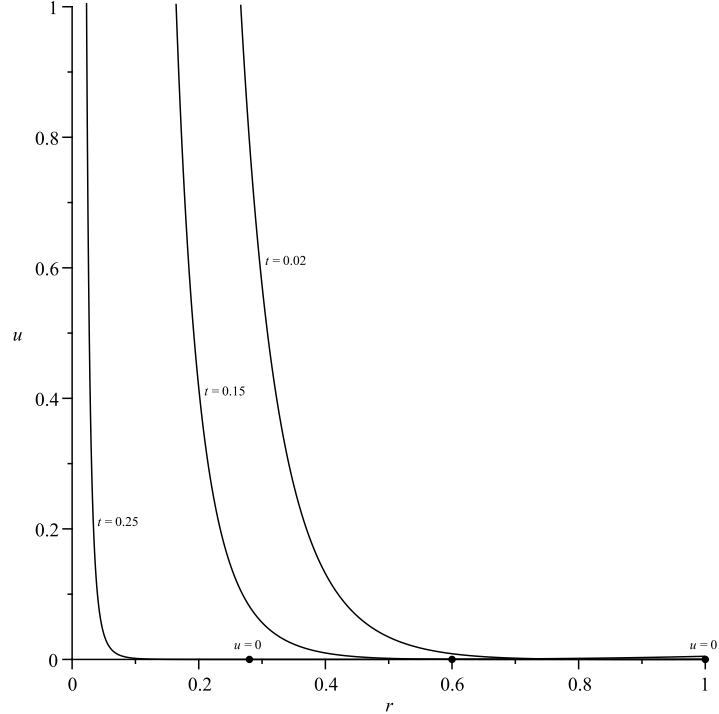}
\end{center}
\caption{Graph of solution \eqref{usol2'} with $n=6~(q+1=1/2).$}\label{eq-sol2}
\end{figure}

At all space-time points away from $r=0$,
this modified similarity solution \eqref{usol2'} is continuous
and has continuous partial derivatives of order $n-3 \geq 2$,
so thus it satisfies the heat equation \eqref{heateq} in the classical sense
(i.e. $u$ is at least $C^1$ in $t$ and $C^2$ in $r$)
on the spatial domain $r>0$.
At $r=0$, $u$ remains singular for all $t>0$,
but since $r^{n-1} u= \beta r t^{n-2} + O(r^2)$ is non-singular
as $r\rightarrow 0$,
$u$ is $n$-dimensionally integrable on the whole spatial domain $r\geq0$ (see figure \ref{eq-sol2}).

A physical interpretation of similarity solution \eqref{usol2'}
can be seen by considering the radial heat flux equation
\EQ
\dfrac{dH}{dt} = S+F
\label{Hprimeeq}
\endEQ
satisfied by the radial heat integral
\EQ
H=\int_0^\infty u r^{n-1} dr
\label{heatintegral}
\endEQ
which is a measure of the total amount of heat in $u(t,r)$,
where
\EQ
F= \lim_{r\rightarrow 0} (-r^{n-1} u_r)
\label{heatflux}
\endEQ
defines the outward radial heat flux at the origin,
and
\EQ
S= k\int_0^\infty u^{q+1} r^{n-1} dr
\label{heatsource}
\endEQ
gives the net amount of heating or cooling produced by
the nonlinear source/sink term in the heat equation \eqref{heateq}.
For the solution \eqref{usol2'} these quantities are given by
\begin{align}
\int_0^{\alpha\sqrt{t}} \beta \left ( t-(r/\alpha)^2 \right)^{n-2} r dr
= \frac{n-4}{n-1}\beta t^{n-1} = H >0,
\label{usol2'H}\\
\lim_{r\rightarrow 0} \left( -(2-n)\beta \left ( t-(r/\alpha)^2 \right)^{n-2} +O(r) \right)
= (n-2)\beta t^{n-2} = F > 0,
\label{usol2'F}\\
\int_0^{\alpha\sqrt{t}} \beta^{(n-4)/(n-2)} \left ( t-(r/\alpha)^2 \right)^{n-4} r^3 dr
= -2\beta t^{n-2} = S< 0.
\label{usol2'S}
\end{align}
Hence $u$ has a positive amount of heat \eqref{usol2'H}
that increases with $t$ due to an increasing, positive
radial outward heat flux at $r=0$.
Since $n\geq 5$,
this flux \eqref{usol2'F} is greater than the net cooling \eqref{usol2'S}
caused by the nonlinear sink term.
Therefore this similarity solution \eqref{usol2'}
physically describes the dispersion of heat
produced by an outward radial heat flux at the origin in $\mathbb{R}^n$,
with some heat absorbed by a nonlinear heat sink proportional to
$u^{(n-4)/(n-2)}$ at all points in $\mathbb{R}^n$.
The dispersion has the behaviour of a radial temperature front at
$r=\alpha\sqrt{t}$ that moves outward with speed $dr/dt = \alpha/(2\sqrt{t})$
for all $t>0$.

Interestingly,
for long times, $u$ increases without bound in $t$
at any spatial point $r>0$ inside the temperature front.
This non-dispersive temporal behaviour is related to $u$ having
an infinite ``energy'' \eqref{energy}, i.e. $E=+\infty$,
so that the flux equation \eqref{energyflux} is not well defined,
allowing $u\rightarrow \infty$ as $t\rightarrow \infty$
despite the ``energy'' integral being formally positive due to $k<0$.

\subsection{ Behavior of solutions for $n\neq 2,3,$ etc. }
\label{nonintegercase}
\noindent\newline\indent
The heat equation \eqref{heateq} can be written in a different form
\EQ\label{2dimheateq}
u_{t} = u_{rr} +(1-\nu)r^{-1}u_{r} + k u^{q+1},
\quad
\nu=\const,\quad
k=\const
\endEQ
in terms of a parameter $\nu=2-n$
which applies to non-integer values of $n$.
This equation \eqref{2dimheateq} describes radial heat flow in $\mathbb{R}^2$
with an extra source/sink term given by $\nu u_r/r$ \cite{Rubinstein-book}.
To interpret this term physically,
we consider the $2$-dimensional radial heat integral
\EQ
H=\int_0^\infty u r dr
\label{2dimheatintegral}
\endEQ
satisfying the radial flux equation
\EQ
\dfrac{dH}{dt} = S+F +\nu\lim_{r\rightarrow 0} u
\label{2dimHprimeeq}
\endEQ
where
\EQ
F= \lim_{r\rightarrow 0} (-r u_r)
\label{2dimheatflux}
\endEQ
defines the outward radial heat flux at the origin,
and
\EQ
S= k\int_0^\infty u^{q+1} r dr
\label{2dimheatsource}
\endEQ
gives the net amount of heating or cooling caused by
the nonlinear source/sink term in the heat equation \eqref{2dimheateq}.
The flux equation \eqref{2dimHprimeeq} shows that,
for the respective cases $\nu >0$ or $\nu <0$,
the term $-\nu u_r/r$ has the interpretation of a heating or cooling point-source
at the origin in $\mathbb{R}^2$.
Thus, solutions $u(t,r)$ will physically describe radial heat flow
arising from a point source in a thin layer,
with heat also produced or absorbed at all points in the layer
due to a nonlinear source/sink term $k u^{q+1}$.

Note that the $2$-dimensional heat equation \eqref{2dimheateq}
retains the form of a gradient flow \eqref{gradientflow} with $n=2-\nu$.

Consider similarity solution \eqref{usol2} with $n=2-\nu$:
\EQ
u(t,r) =\beta r^{\nu} \left( \pm(\alpha(t+c) +r^2) \right)^{-\nu},
\quad q=-2/\nu, \quad \nu\neq 0,-1,-2
\label{2Dusol2}
\endEQ
where
\EQ
\alpha = 2(\nu+2) ,\quad
\beta =(-4\nu(\nu+1)/k)^{\nu/2} ,
\endEQ
which requires $k>0$ if $0>\nu>-1$ and $k<0$ if $\nu<-1$ or $\nu>0$.
The behaviour of this solution depends essentially on the separate signs of
$\alpha$ and $\nu$.

For $\alpha <0$, \eqref{2Dusol2} can be modified similarly to \eqref{usol2'}
by putting a cutoff on $u$ at the space-time parabola where $u=0$.
Then
\EQ
u(t,r) =\begin{cases}
0, & r\geq \sqrt{t/|\alpha|} \\
\beta r^{\nu} \left( |\alpha|t-r^2 \right)^{-\nu}, & 0\leq r\leq \sqrt{t/|\alpha|}
\end{cases}
\label{2Dusol2'}
\endEQ
gives a classical solution of the heat equation \eqref{2dimheateq}
on the spatial domain $r>0$
(i.e. $u$ belongs to $C^1(\mathbb{R}^+)$ in $t$ and $C^2(\mathbb{R}^+)$ in $r$),
with $\nu<-2$ and $k<0$.
However,
at $r=0$, $u$ is singular such that, for all $t>0$,
$H=+\infty$
and
$\displaystyle F +\nu\lim_{r\rightarrow 0} u =
2(\nu/\gamma) t^{-\nu-1} \lim_{r\rightarrow 0} r^{\nu+1} = +\infty$
since $\nu+1<-1$.
The modified similarity solution \eqref{2Dusol2'} therefore has
the physical interpretation of heat dispersion produced by
an infinite net outward radial heat flux at the origin in $\mathbb{R}^2$,
with a radial temperature front located at $r=\sqrt{t/|\alpha|}$ for $t>0$.

For $\alpha >0$, \eqref{2Dusol2} is smooth and positive
on the spatial domain $r>0$.
It thus gives a $C^\infty(\mathbb{R}\times\mathbb{R}^+/\{0\})$ solution
of the heat equation \eqref{2dimheateq}
with the asymptotic behaviour
\EQ
u(t,r) =\beta r^{\nu} (\alpha t +r^2)^{-\nu}
= \begin{cases}
r^{\nu} (\beta\alpha^{-\nu} t^{-\nu} +O(r)), & r\rightarrow 0 \\
r^{-\nu} (\beta +O(1/r)), & r\rightarrow \infty \\
\end{cases}
\label{2Dusol2''}
\endEQ
(with $c=0$) where $\nu>-2$.
In the least interesting case when $0>\nu>-2$,
$u$ is unbounded for large $r$,
whereby it has an infinite amount of heat $H=+\infty$
and ``energy'' $E=+\infty$.
In contrast, when $\nu >0$,
$u$ decays to $0$ for large $r$ and vanishes at $r=0$ for $t>0$.
As a result, in this case the ``energy'' of $u$ for $t>0$ is given by
\EQ
E = \frac{1}{2^{3\nu/2}3} \nu (\nu^2+2\nu+6) (\nu+2)^{(1-\nu)/2}
\frac{\Gamma(1+\nu/2)\Gamma(3\nu/2)}{\Gamma(2+2\nu)} \beta t^{-(1+3\nu)/2}
\endEQ
which is finite and positive due to $k<0$ for $\nu>0$.
In particular, $E$ decreases to $0$ as $t\rightarrow \infty$
in accordance with the flux equation \eqref{energyflux}.
Moreover, provided $\nu>2$,
$u$ has a finite amount of heat
\EQ
H = \frac{4\sqrt{\pi}}{2^{3\nu/2}3} (\nu+2)^{(1-\nu)/2} (\nu-2)^{-1}
\frac{\Gamma(1+\nu/2)}{\Gamma((1+\nu)/2)} \beta t^{1-\nu/2}
\endEQ
which decreases to $0$ as $t\rightarrow \infty$.
For $t>0$ the corresponding heat flux quantities in this case are given by
\begin{align}
S = -\frac{\nu-2}{2} H/t < 0,
\\
F=0 ,\quad
\nu\lim_{r\rightarrow 0}u =0 .
\end{align}

\begin{figure}[h]
\begin{center}
\includegraphics[scale=0.5]{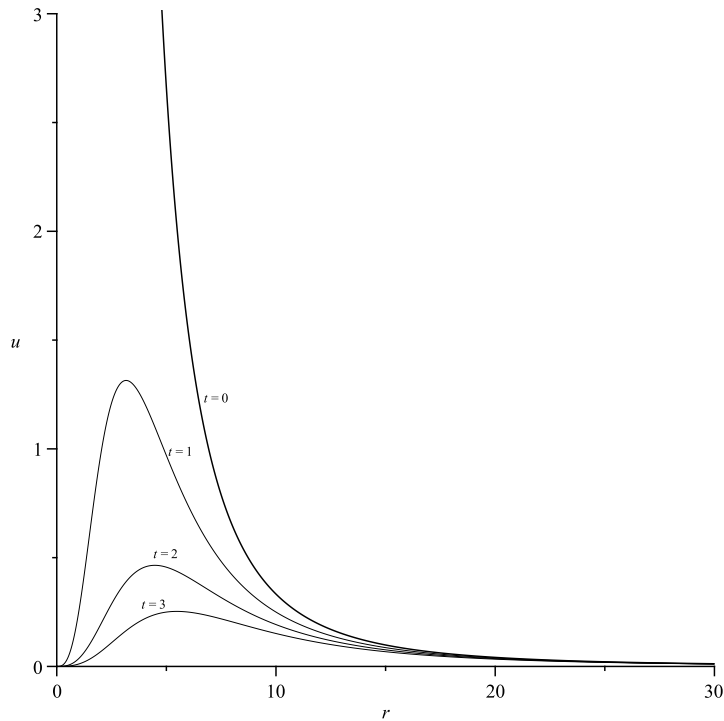}
\end{center}
\caption{Graph of solution \eqref{2Dusol2''} with $\nu=3$.}\label{eq-2Dsol2}
\end{figure}

Therefore, in the most physically interesting case $\nu>2$,
the similarity solution \eqref{2Dusol2''} describes
an initial monopole-like circular heat distribution in $\mathbb{R}^2$ (see figure \ref{eq-2Dsol2})
producing a smooth positive dispersive radial heat flow
with some heat absorbed by a nonlinear heat sink proportional to
$u^{1-2/\nu}$ at all points in $\mathbb{R}^2$.
In particular,
this solution has no point-source or heat flux at the origin for all times $t>0$.

Similarity solutions \eqref{usol5} and \eqref{usol6} have
$q+1=3$, $\nu=-1/2$, and $k<0$.
In both solutions, for all $t$,
$u$ is singular at $r=0$ and has slow radial decay such that $u=O(1/r)$ for large $r$.
Consequently, the amount of heat in $u$ is $H=\pm\infty$.
In particular,
at $t=0$, $u$ is given by a $1/r$ heat distribution for $r\geq 0$,
while for $t\rightarrow \infty$ at any spatial point $r>0$,
$u=r^{-1}( \pm(1/\sqrt{-2k}) +O(1/t) )$
approaches a constant multiple of the same $1/r$ heat distribution.
This non-dispersive temporal behaviour occurs because
$u$ has ``energy'' $E=+\infty$, which does not decrease with $t$.

\begin{figure}[h]
\begin{center}
\includegraphics[scale=0.5]{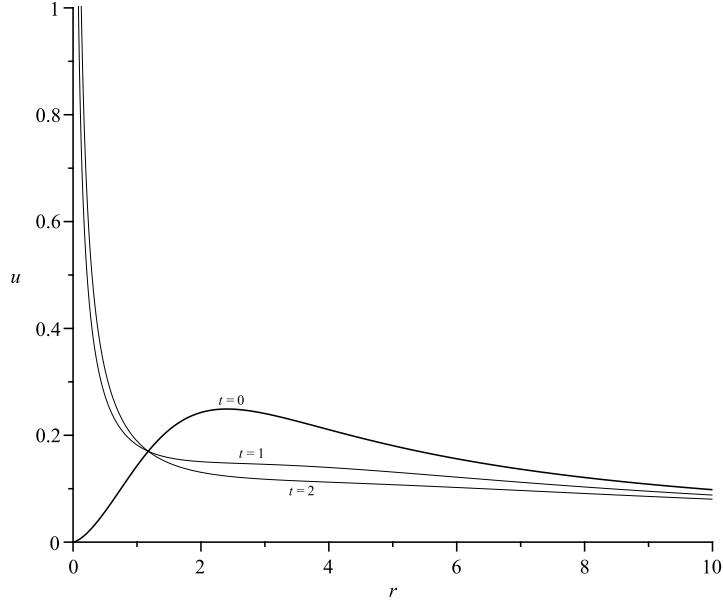}
\end{center}
\caption{Graph of solution \eqref{usol4}.}\label{eq-sol4}
\end{figure}

Non-similarity solution \eqref{usol4} also has $q+1=3$, $\nu=-1/2$, and $k<0$.
Compared with the previous similarity solutions \eqref{usol5} and \eqref{usol6},
it exhibits the same long-time behaviour
$u=r^{-1}( \pm(1/\sqrt{-2k}) +O(1/t) )$
as $t\rightarrow \infty$ at any spatial point $r>0$.
It also exhibits the same radial decay $u=O(1/r)$ as $r\rightarrow\infty$
for all $t$, whereby $H=\pm\infty$ and $E=+\infty$.
However,
near $r=0$, \eqref{usol4} is less singular than \eqref{usol5} and \eqref{usol6},
such that $u=r^{-1/2}( \pm(15/c\sqrt{-2k})t +O(r) )$.
As a consequence the heat flux quantities for \eqref{usol4}
have the properties $F + \nu\lim_{r\rightarrow 0}u =0$ and $0<-S <\infty$.
Moreover, at $t=0$, \eqref{usol4} reduces to the heat distribution
$u=\pm(5/\sqrt{-2k})(r+c/r^{3/2})^{-1}$
which is $C^1$ for $r\geq0$ if $c>0$
and which vanishes at both $r=0$ and $r=\infty$ (see figure \ref{eq-sol4}).
This non-similarity solution \eqref{usol4}
therefore has the physical interpretation of
an initial circularly peaked heat distribution
which vanishes both at the origin and spatial infinity in $\mathbb{R}^2$,
producing a radial heat flow that is singular at the origin
and has a $1/r$ decay for large radius,
with a finite amount of heat absorbed by a nonlinear heat sink proportional to
$u^3$ at all points in $\mathbb{R}^2$.
The singularity in the heat distribution $u$ at the origin corresponds to
an infinite cooling point-source plus an infinite outward heat flux,
whose net contribution to the heat flow for all times $t> 0$ is zero.

\begin{figure}[h]
\begin{center}
\includegraphics[scale=0.5]{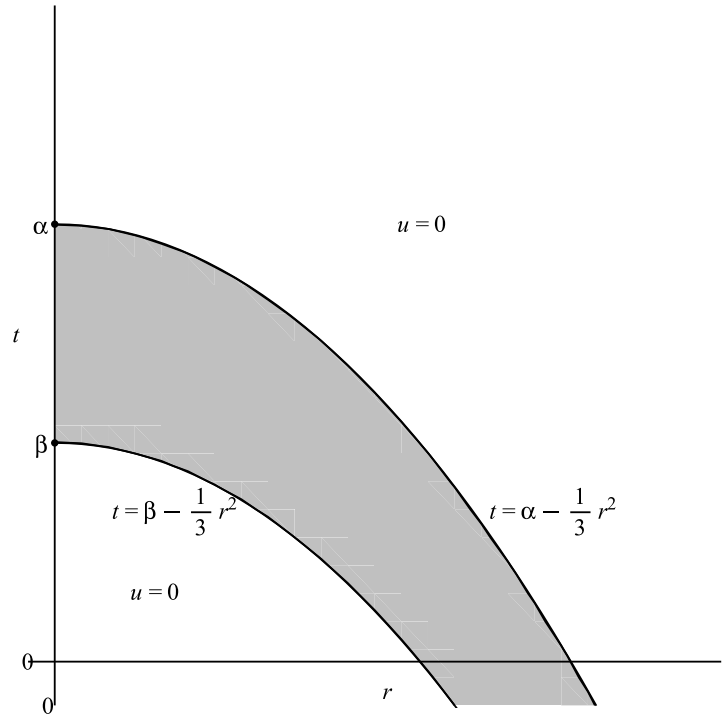}
\end{center}
\caption{Space-time graph of solution \eqref{2Dnonsimsol1}.}\label{eq-2Dnonsimsol1-space-time}
\end{figure}

\begin{figure}[h]
\begin{center}
\includegraphics[scale=0.5]{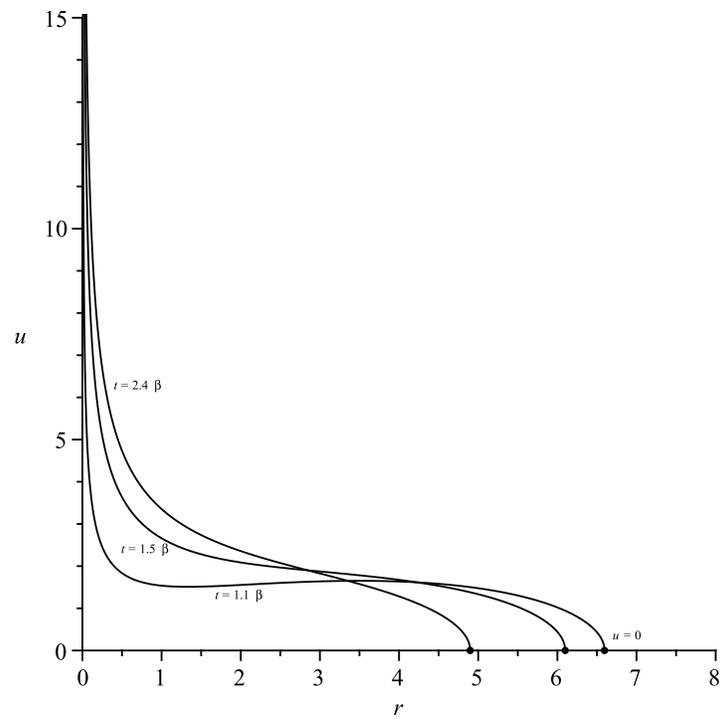}
\end{center}
\caption{Graph of solution \eqref{2Dnonsimsol1} with $\beta=5,~ \alpha=20$.}\label{eq-2Dnonsimsol1-beta-t-alpha}
\end{figure}

\begin{figure}[h]
\begin{center}
\includegraphics[scale=0.5]{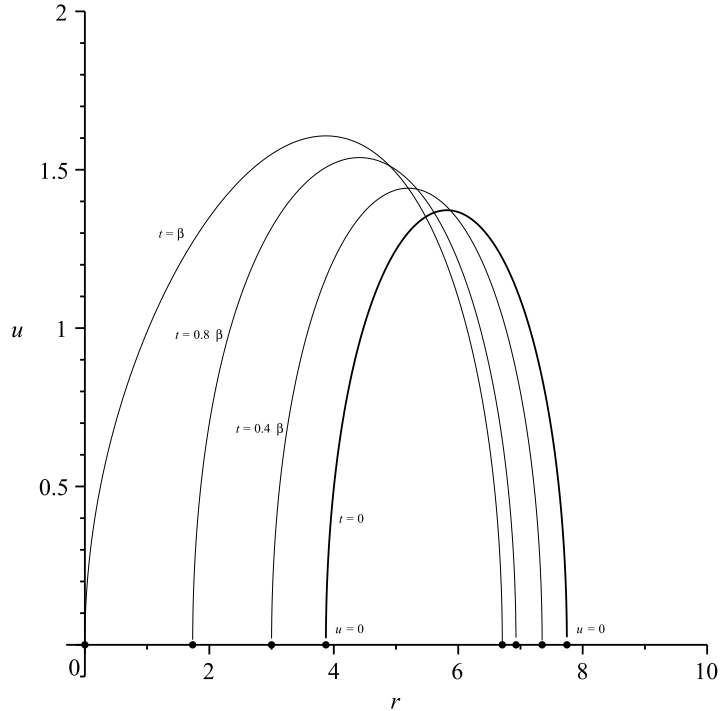}
\end{center}
\caption{Graph of solution \eqref{2Dnonsimsol1} with $\beta=5$.}\label{eq-2Dnonsimsol1-t-less-beta}
\end{figure}

Finally, non-similarity solution \eqref{nonsimsol1} has
$q+1=-3$, $\nu=-1/2$, and $k>0$.
In this solution, for all $t$,
$u$ is unbounded as $r\rightarrow\infty$ and has a cusp
at the space-time parabolas $t+\tilde c=-r^2/3$ and $t+\tilde c=-(r^2+1/c)/3$
where $u_r$ blows up and $u$ vanishes. 
We can modify the solution by putting a cutoff on $u$ at both parabolas so that
\EQ
u(t,r) = \begin{cases}
0 , & t-\beta\leq -r^2/3 \mathtext{ or } t-\alpha\geq -r^2/3 \\
r^{-1/2} \sqrt{\gamma (3(\alpha-t)-r^2) (3(t-\beta)+r^2)},
& t-\beta\geq -r^2/3 \mathtext{ and } t-\alpha\leq -r^2/3 \\
\end{cases}
\label{2Dnonsimsol1}
\endEQ
(see figure \ref{eq-2Dnonsimsol1-space-time})
where
\EQ
\gamma = \frac{\sqrt{k}}{3(\alpha-\beta)}
\endEQ
and $\alpha=-\tilde c>\beta=-(\tilde c +1/3c)>0$ with $c>0$.
Then the modified solution \eqref{2Dnonsimsol1} is $C^\infty$ in $r$ and $t$
at all space-time points other than $r=0$,
$\alpha-t=r^2/3$ and $\beta-t=r^2/3$.
Near $r=0$, for $\beta<t<\alpha$,
$u=r^{-1/2} \left( 3\sqrt{\gamma (\alpha-t)(t-\beta)} +O(r)\right)$
is singular (see figure \ref{eq-2Dnonsimsol1-beta-t-alpha}),
while near the inner and outer cusps,
$u=(\sqrt{3(\alpha-t)}-r)^{1/2} \left( \sqrt{6\gamma  (\alpha-\beta)} +O(\sqrt{3(\alpha-t)}-r) \right)$ for $t<\alpha$
and
$u=(r-\sqrt{3(\beta-t)})^{1/2} \left( \sqrt{6\gamma  (\alpha-\beta)} +O(r-\sqrt{3(\beta-t)}) \right)$ for $t<\beta$
each have square-root behaviour in $r$ (see figure \ref{eq-2Dnonsimsol1-t-less-beta}).
From these properties, \eqref{2Dnonsimsol1} can be checked to satisfy
the heat equation \eqref{2dimheateq} in a weak sense
(i.e. $u$ is only $C^0$ in $t$ and $r$) on the spatial domain $r\geq 0$.
In particular, $u$ has a finite, non-negative amount of heat
\begin{align}
H=\begin{cases}
0,
& t>\alpha \\
\sqrt{\gamma } \displaystyle\int_0^{\sqrt{3(\alpha-t)}} (3(\alpha-t)-r^2)^{1/2} (3(t-\beta)+r^2)^{1/2} r^{1/2} dr
\\\quad\displaystyle
= \frac{\pi}{16} k^{1/4} (3(\alpha-\beta))^{3/2} (3(\alpha-t))^{-1/4} {}_2F{}_1(1/4,3/2;3;(\alpha-\beta)/(\alpha-t)) ,
& \beta<t<\alpha
\\
\sqrt{\gamma } \displaystyle\int_{\sqrt{3(\beta-t)}}^{\sqrt{3(\alpha-t)}} (3(\alpha-t)-r^2)^{1/2} (3(t-\beta)+r^2)^{1/2} r^{1/2} dr
\\\quad\displaystyle
=\frac{2}{5}\sqrt{\frac{2}{\pi}} \Gamma(3/4)^2 k^{1/4} (3(\alpha-t))^{5/4} {}_2F{}_1(-1/2,3/2;9/4;(\alpha-t)/(\alpha-\beta)) ,
& t<\beta
\\
\end{cases}
\label{2Dnonsimsol1H}
\end{align}
so thus $u$ is in $L^1(\mathbb{R}^+)$ for all $t\geq 0$.
The corresponding heat flux equation is given by
\EQ
\dfrac{dH}{dt} = \begin{cases}
0,
& t>\alpha \\\displaystyle
F-\lim_{r\rightarrow 0} u/2 -F_{\rm outer} +S
\\\quad\displaystyle
= -\frac{3\pi}{8} k^{1/4} (3(\alpha-\beta))^{1/2} (3(\alpha-t))^{-1/4}
\big( {}_2F{}_1(1/4,1/2;2;(\alpha-\beta)/(\alpha-t))
\\\qquad\qquad
- {}_2F{}_1(1/4,3/2;2;(\alpha-\beta)/(\alpha-t)) \big) ,
& \beta<t<\alpha
\\
-F_{\rm inner} +F_{\rm outer} +S
\\\quad\displaystyle
=-\frac{3}{5}\sqrt{\frac{2}{\pi}} \Gamma(3/4)^2 k^{1/4} (3(\alpha-t))^{1/4}
\big( \frac{5}{2}\ {}_2F{}_1(-1/2,3/2;5/4;(\alpha-t)/(\alpha-\beta))
\\\qquad\qquad\displaystyle
- \frac{\alpha-t}{\alpha-\beta}\ {}_2F{}_1(1/2,3/2;9/4;(\alpha-t)/(\alpha-\beta)) \big) ,
& t<\beta
\\
\end{cases}
\label{2Dnonsimsol1fluxeq}
\endEQ
where
\begin{align}
&
F-\lim_{r\rightarrow 0} u/2
= \lim_{r\rightarrow 0}(ru_r-u/2) = 0,
\\
&
F_{\rm inner} = \lim_{r\rightarrow \sqrt{3(\beta-t)}} ru_r = +\infty,
\\
&
F_{\rm outer} = \lim_{r\rightarrow \sqrt{3(\alpha-t)}} ru_r = -\infty,
\\
&
S = \begin{cases}
k \gamma^{-3/2} \displaystyle\int_0^{\sqrt{3(\alpha-t)}}
(3(\alpha-t)-r^2)^{-3/2} (3(t-\beta)+r^2)^{-3/2} r^{-1/2} dr,
& \beta<t<\alpha \\
k\gamma^{-3/2} \displaystyle\int_{\sqrt{3(\beta-t)}}^{\sqrt{3(\alpha-t)}}
(3(\alpha-t)-r^2)^{-3/2} (3(t-\beta)+r^2)^{-3/2} r^{-1/2} dr,
& t<\beta . \\
\end{cases}
\end{align}
Here the quantities $F_{\rm inner}$ and $F_{\rm outer}$ are defined to be
inward/outward heat fluxes arising from the inner and outer cusps, respectively.
These fluxes act as cooling sources that cancel the singular contributions
coming from the endpoints in the nonlinear heating source $S$
in the flux equation \eqref{2Dnonsimsol1fluxeq}.
Their net effect gives a finite cooling rate, $dH/dt<0$.
Hence $H$ decreases from its initial value $H(0)\neq 0$ at $t=0$
to $H(\alpha)=0$ at $t=\alpha$.

\clearpage

As a result, non-similarity solution \eqref{nonsimsol1}
physically describes an initial ring-shaped heat distribution
in $\mathbb{R}^2$ with radial temperature fronts at
$r=\sqrt{3(\beta-t)}$ and $r=\sqrt{3(\alpha-t)}$ where $0<\beta<\alpha$.
These fronts behave as heat flux sinks and move radially toward the origin
at speeds $dr/dt= -\sqrt{3/(\beta-t)}$ and $dr/dt= -\sqrt{3/(\alpha-t)}$.
At time $t=\beta$ the inner front reaches the origin and produces
an infinite cooling point-source plus a compensating infinite outward heat flux
that both persist until time $t=\alpha$
when the outer front reaches the origin.
For all times $0\leq t\leq\alpha$,
the heat absorbed at the temperature fronts exceeds the heat produced by
the nonlinear source term $u^{-3}$
at all points in $\mathbb{R}^2$ inside the fronts,
whereby the total amount of heat decreases to zero in a finite time $t=\alpha$,
with the heat distribution being given by $u=0$ at all points $r\geq 0$
for times $t\geq\alpha$.

\section{ Concluding remarks }
\label{remarks}

As main results in this paper,
analytically interesting exact solutions have been obtained for
a multi-dimensional semilinear heat equation \eqref{heatfloweq}
via a separation technique applied to the group foliation equations
associated with the group of scaling symmetries \eqref{scaling}
admitted by this equation.
The solutions consist of explicit similarity solutions as well as
other explicit solutions of a more general (non-similarity) form.

In general our method provides a highly effective alternative
to standard similarity reduction for finding exact solutions
to nonlinear PDEs with a group of scaling symmetries.
Firstly,
this method is an algorithmic refinement of the basic approach
developed for the semilinear wave equation in \cite{AncoLiu},
which leads to systematic reductions of the group foliation equations
into overdetermined systems of ODEs that can be derived and solved
by means of computer algebra.
In particular, these ODE systems are tractable to solve
using the computer algebra package {\sc Crack} \cite{crack}.

Secondly,
the method is able to yield exact similarity solutions in an explicit form,
whereas standard similarity reduction only gives an ODE that
still has to be solved to find solutions explicitly
and in general this step can be quite difficult.
Indeed, the resulting similarity ODE \eqref{similarityeq}
for the semilinear heat equation \eqref{heatfloweq}
cannot be solved by standard integration techniques
such as symmetry reduction or integrating factors.

Thirdly,
because the group foliation equations contain
all solutions of the given nonlinear PDE,
our method can yield non-similarity solutions that are also not invariant
under any other (non-scaling) point symmetries admitted by the nonlinear PDE.

We can apply the same method more generally to
nonlinear PDEs without scaling symmetries
by utilizing the group foliation equations associated with
any admitted one-dimensional group of point symmetries of
the given nonlinear PDE
and by adapting the separation technique to the specific form of
the non-derivative terms that appear in the given group foliation equations.
The algorithmic aspects of these steps will be the same as
in the similarity case we have presented in this paper,
since every one-dimensional group of point symmetries can be equivalently
expressed as a group of scalings under an appropriate
change of independent and dependent variables
(i.e. by an invertible point transformation).

For future work, we plan to present a full comparison between
the present group foliation method and standard symmetry reduction
as applied to many typical linear and nonlinear PDEs of interest, e.g.
linear heat and wave equations;
semilinear diffusion and telegraph equations;
integrable semilinear evolution equations
such as the Korteweg de Vries and Boussinesq equations.

\section{Acknowledgement}
S. Anco and T. Wolf are each supported by an NSERC research grant.

S. Ali thanks the Mathematics Department of Brock University for support
during the period of a research visit when this paper was written.

Computations were partly performed on computers of the Sharcnet consortium (www.sharcnet.ca).

The referees are thanked for valuable comments which have improved this paper.

\end{document}